\documentclass[12pt, reqno]{amsart}

\author[M.~Caprio and S.~Mukherjee]{Michele Caprio and Sayan Mukherjee}
\address{Department of Computer  Science, University of Manchester, Oxford Road, Manchester M13 9PL, United Kingdom}
\email{michele.caprio@manchester.ac.uk}
\urladdr{\url{https://michelecaprio.wixsite.com/caprio}}   
\address{Max Planck Institute for Mathematics in the Sciences, Inselstraße 22, Leipzig 04103, Germany}
\email{sayan.mukherjee@mis.mpg.de}
\urladdr{\url{https://sayanmuk.github.io/}}


\title{Finite Admixture Models: a Bridge with Stochastic Geometry and Choquet Theory}

\usepackage[T1]{fontenc}
\usepackage{amsmath}
\usepackage{amssymb}
\usepackage{amsthm}
\usepackage[left=2.5cm, right=2.5cm, top=3cm]{geometry}
\usepackage{hyperref}
\usepackage{algorithm}
\usepackage{algorithmicx}
\usepackage{algpseudocode}
\usepackage{bbm}
\usepackage{tikz}
\usepackage{caption}
\usepackage{subcaption}
\usepackage{fancyhdr}
\usepackage[inline]{enumitem}
\usepackage{comment}
\usepackage{nicefrac}
\usepackage{bm}
\usepackage{mathrsfs}
\usepackage{graphicx}
\usepackage[utf8]{inputenc}
\usepackage{cancel}
\usepackage{mathtools}

\newcommand{\vertiii}[1]{{\left\vert\kern-0.25ex\left\vert\kern-0.25ex\left\vert #1 
    \right\vert\kern-0.25ex\right\vert\kern-0.25ex\right\vert}}
		
\AtBeginDocument{%
   \def\MR#1{}
}

\algdef{SE}[DOWHILE]{Do}{doWhile}{\algorithmicdo}[1]{\algorithmicwhile\ #1}%

\newtheorem{thm}{Theorem}[section]

\theoremstyle{definition} 
\newtheorem{defi}[thm]{Definition}
\let\olddefi\defi
\renewcommand{\defi}{\olddefi\normalfont}
\newtheorem{rmk}[thm]{Remark}
\let\oldrmk\rmk
\renewcommand{\rmk}{\oldrmk\normalfont}

\DeclareMathOperator*{\argmax}{arg\,max}

\newtheorem{theorem}{Theorem}

\newtheorem{proposition}[theorem]{Proposition}
\newtheorem{corollary}{Corollary}[theorem]
\newtheorem{claim}[theorem]{Claim}
\newtheorem{definition}[theorem]{Definition}
\newtheorem{remark}{Remark}
\newtheorem{example}[theorem]{Example}

\hypersetup{
    pdftitle={prove},
    pdfauthor={autori},
    pdfmenubar=false,
    pdffitwindow=true,
    pdfstartview=FitH,
    colorlinks=true,
    linkcolor=blue,
    citecolor=green,
    urlcolor=cyan
}

\providecommand{\MR}[1]{}

\providecommand{\MR}{\relax\ifhmode\unskip\space\fi MR }

\providecommand{\href}[2]{#2}

\newcommand\xqed[1]{%
  \leavevmode\unskip\penalty9999 \hbox{}\nobreak\hfill
  \quad\hbox{#1}}
\newcommand\demo{\xqed{$\triangle$}}

\begin{document}

\begin{abstract}
Given a finite admixture model whose components and weights are unknown, let the number of identifiable components be a function of the amount of data sampled from a known distribution on the unit simplex. We use techniques from stochastic convex geometry to find the growth rate of its expected value. In addition, when the components are known but the weights are not, we provide an application of the classic Glivenko-Cantelli's theorem that allows us to retrieve the Choquet measure supported on the identifiable admixture components. In turn, this gives us the identifiable admixture weights. Finally, we propose a novel algorithm that estimates the model capturing the complexity of the data using only the strictly necessary number of components.
\end{abstract}

\maketitle
\thispagestyle{empty}

\section{Introduction}\label{intro}
Finite mixture models go back at least to \cite{pearson2,pearson1} and have served as a workhorse in stochastic modeling \cite{everitt,lindsay,mengersen}. Applications include clustering \cite{mclachlan1}, hierarchical or latent space models \cite{masyn}, and semiparametric models \cite{mcnicholas} where a mixture of simple distributions is used to model data that is putatively generated from a complex distribution. In finite mixture models, the mixing distribution is over a finite number of components; there are also many examples of infinite mixture models in the Bayesian nonparametrics literature \cite{antoniak1974,west94}.

We consider a finite mixture of multinomials. We start with the basic multinomial model where an observation $X$ takes on $J$ possible values $\{1,\ldots,J\}$, and
$X \sim \mbox{Mult}(\pi)$, with $\pi \equiv (\pi_1,\ldots,\pi_J)^\top$ where $\pi_j=\mathbb{P}(X= j)$, with $\pi_j \geq 0$ for all $j$ and $\sum_{j=1}^J \pi_j=1$. Notice that here and in the remainder of the paper, we adopt the notational convention that, when not specified, the number of trials for a Multinomial distribution is equal to $1$. A mixture of $L$ multinomials can be specified as follows
\begin{equation}\label{fam_first}
    X \sim \mbox{Mult}(\pi), \quad \pi = \sum_{\ell=1}^L \phi_{\ell} f_\ell,
\end{equation}
where probability vector $\phi \equiv (\phi_{1},\ldots,\phi_{L})^\top$ assigns the probability of the  observation $X$ coming from the $\ell$-th mixture component with multinomial parameter 
$$f_\ell=(f_{\ell,1},\ldots,f_{\ell,J})^\top.$$
We have that $\sum_{j=1}^J f_{\ell,j} =1$ with $f_{\ell,j} \geq 0$, and $\sum_{\ell=1}^L \phi_{\ell} =1$ with $\phi_{\ell} \geq 0$, for all $\ell \in \{1,\ldots,L\}$. An important point throughout the paper is that $\pi$ belongs to the convex hull of probability vectors  $\{f_1,\ldots,f_L\}$. The convex hull of $\{f_1,\ldots,f_L\}$ is a function of the {\em identifiable elements} of $\{f_1,\ldots,f_L\}$, that is, {\em those elements that cannot be written as a convex combination of the other $f_\ell$'s}. In other words, we call identifiable the {\em extremal elements} of the convex hull. Hence, understanding the identifiable elements provides information about the  key model parameter $\pi$.
 

The finite mixture model we stated is an example of a finite admixture model; the most popular finite admixture model is the latent Dirichlet allocation (LDA) model \cite{blei, Pritchard945}. A classic application of an admixture model is a generative process for documents. Consider a document as a collection of words; LDA posits  that each document is a mixture of a small number of topics, and that these latter can be modeled by a multinomial distribution on the presence of a word in the topic. The hierarchical Dirichlet process \cite{teh}, and generalizations thereof, may be considered as the natural nonparametric counterpart of the LDA model.

The $f_\ell$'s and $\pi$ in \eqref{fam_first} are all elements of $\Delta^{J-1}$, the unit simplex on $\mathbb{R}^J$. Again, $\pi$ belongs to the convex hull of  $\{f_\ell\}_{\ell=1}^L$, or $\pi \in \text{Conv}(f_1,\ldots,f_L)$. Hence, an element of a convex hull in the Euclidean unit simplex represents (the distribution of) a finite admixture model. 

Notice that the number of extremal elements of $\text{Conv}(f_1,\ldots,f_L)$, which we denote as $M$, will probably be less than $L$ because some of the components $f_\ell$ are likely to be a convex combination of the others. A key concept in this paper is what we call the \textit{richest cheap model} representing $\pi$, that is, the finite admixture model representing $\pi$ whose components are $\{f_k\}_{k \in \mathcal{I}}$ such that $f_k \not\in \text{Conv}(f_{\mathcal{I}\setminus \{k\}})$, for all $k \in \mathcal{I}$, $\mathcal{I} \subset \{1,\ldots,L\}$, and $\#\mathcal{I}=M$, where $\#$ denotes the cardinality operator. These conditions tell us that the $M$ components of the richest cheap model are a subset of $\{f_1,\ldots,f_L\}$ and cannot be written as a convex combination of one another. By assuming -- without loss of generality -- that the identifiable elements in $\{f_1,\ldots,f_L\}$ are the first $M$ ones, for the richest cheap model we can rewrite $\pi$ in \eqref{fam_first} as
\begin{equation}\label{rcm}
    \pi=\sum_{\ell=1}^M \varphi_{\ell} f_\ell,
\end{equation}
where we denote by $\varphi \equiv (\varphi_{1},\ldots,\varphi_{M})^\top$ the probability vector that assigns the probability of the observation $X$ coming from the $\ell$-th identifiable component with multinomial parameter $f_\ell=(f_{\ell,1},\ldots,f_{\ell,J})^\top$, $\ell\in\{1,\ldots,M\}$. That is, the $\varphi_{\ell}$'s are the identifiable admixture weights. Of course the $\varphi_{\ell}$'s are such that $\sum_{\ell=1}^M \varphi_{\ell}=1$, and $\varphi_{\ell} \geq 0$, for all $\ell \in \{1,\ldots,L\}$. As we can see, the richest cheap model captures the underlying complexity associated with the data at hand, using only the strictly necessary number of components.

Rather than developing new tools for working with or applying finite admixture models, the main goal of this paper is to establish connections between finite admixture models, Choquet theory, and stochastic convex geometry. This paper is an extension of \cite[Chapter 2]{thesis}.


\subsection{Choquet theory}
Choquet theory, named after mathematician Gustave Choquet, is an area of functional and convex analyses concerned with measures which have support on the extreme points of a convex set \cite{phelps1}. Its fundamental tenet is that we can represent every element in a convex set $C$ via a weighted average of the extremal elements of the set. Here, weighted average is to be understood as a generalization of the usual notion of convex combination to an integral taken over the set $E$ of extreme points of $C$. The formal, central result to Choquet theory is the following.

\begin{theorem}{\textbf{(Choquet, cf. \cite[Theorem, page 14]{phelps1}).}}\label{choq1}
Let $C$ be a metrizable compact convex subset of a locally convex space $V$. Pick any $c \in C$. Then, there exists a probability measure $\nu$ on $C$ which represents $c$ and is supported on $E$, that is,
$$f(c)=\int_E f(e) \nu(\text{d}e),$$
for any affine function $f$ on $C$.
\end{theorem}

Choquet also characterized those compact convex sets $C$ with the property that for every $c \in C$ there is a \textit{unique} probability measure $\nu_c$ supported on $E$ that represents $c$. 


\begin{theorem}{\textbf{(Choquet, cf. \cite[Theorem, page 60]{phelps1}).}}\label{choq2}
Let $C$ be a metrizable closed convex subset of a locally convex space $V$. Then, $C$ is a simplex\footnote{Choquet's result refers to a more general notion than the classical simplex, called {\em Choquet simplex}. Since a classical simplex is a Choquet simplex, and since in this paper we only work with the former, we avoid introducing the notion of Choquet simplex to refrain from unnecessary complications. The interested reader can find it in \cite[Chapter 10, page 52]{phelps}.} if and only if for every $c$ in $C$, there exists a unique measure $\nu_c$ which represents $c$ and is supported on $E$, that is, 
$$f(c)=\int_E f(e) \nu_c(\text{d}e),$$
for any affine function $f$ on $C$.
\end{theorem}

We call $\nu_c$ the \textit{Choquet measure} for $c$. Notice that we differentiate between $\nu$ in Theorem \ref{choq1} and $\nu_c$ in Theorem \ref{choq2} to highlight the fact that the latter uniquely represents $c$. These results entail that studying the extremal elements of a convex set gives us important results concerning the (elements of the) whole set. Choquet theory in the context of finite mixture models has been inspected in \cite{hoff}. There, the author develops an approach that uses Choquet's theorems for inference with the goal of estimating probability measures constrained to lie in a convex set, for example mixture models. The key observation in \cite{hoff} is that  inference over a convex set of measures can be made via unconstrained inference over the set of extreme measures. The main difference between this work and the approach developed in \cite{hoff} is that we consider a convex hull of points in a unit simplex rather than the convex hull of probability measures. Furthermore, our goal is different: we use a result from Choquet theory to retrieve the identifiable weights in the finite admixture model at hand. Notice also that de Finetti's theorem \cite{definetti3,definetti1,definetti2} can be given a geometric interpretation -- inspected in Appendix \ref{infiniteex} -- that is heuristically similar to that of Choquet theory.

\subsection{Stochastic convex geometry}
This paper also establishes a bridge between finite admixture models and stochastic convex geometry that allows to view finite admixture models as well-studied geometric objects. This insight allows to closely relate the number of identifiable admixture components to the number of extremal elements of a convex body. Thereby, it facilitates studying the asymptotic growth rate and the asymptotic distribution of the number of components. 

The geometry of finite mixture models has primarily been studied in two contexts: differential geometry \cite{amari,kass_vos} and convex geometry \cite{lindsay,marriott}. The approach in this paper is based on (stochastic) convex geometry. The first to study the geometry of mixture models was Lindsay \cite{lindsay_classic1,lindsay_classic2}. In the first paper, the author  established the geometric properties of the likelihood set and used these properties to study the uniqueness of the maximum likelihood estimator (MLE) as well as other fundamental properties of the MLE. In the second paper, the author established the results for the nonparametric MLE. Lindsay also wrote a book \cite{lindsay}  
whose focus was the identifiability of the mixture weights, a Carath\'eodory representation theorem for multinomial mixtures, and the asymptotic mixture geometry. In a recent paper \cite{nguyen}, the author studies the asymptotic behavior of the convex polytope representing a finite admixture model. Their results well complement the ones in section \ref{mainresults}, where the focus is the identifiable admixture components, whose geometric representation is given by the extremal elements of the convex polytope representing the finite admixture model. In \cite{marriott}, the author bridges the differential and convex geometric approaches to identify restrictions for which a mixture model can be written as a tractable geometric quantity that can simplify inference problems. This paper is similar in spirit to Lindsay's work, but uses more modern techniques.

\subsection{Main results and structure of the paper}
We provide three main results. The first two, Theorem \ref{growth} and Theorem \ref{extrema2}, state the following. Suppose we do not know what the components and the weights in our admixture model are, and we also do not know the number of components. Then, if we assume that the number of identifiable components $M$ is a function $M(n)$ of the amount $n$ of data that we sample from a known distribution on $\Delta^{J-1}$, we are able to tell the speed at which its expected value grows. The other main result is Theorem \ref{thm1}. It states that if we know the number of identifiable components of the model, but not the components themselves nor the weights, we can apply the classic Glivenko-Cantelli's theorem to retrieve the Choquet measure on the components. In turn, the latter gives us the identifiable admixture weights. 
We also show how looking for the richest cheap model can be seen as an optimization problem, and we propose an algorithm to solve it. 

The paper is organized as follows. In section \ref{mainresults}, we let the number of identifiable admixture components $M$ depend on the sample size $n$, that is, we let $M=M(n)$. 
We state and prove Theorems \ref{growth} and \ref{extrema2}. In addition, in Theorems \ref{clt_2} and \ref{clt_gen} we state a central limit theorem (CLT) for the distribution of the number of identifiable admixture components, and in Theorem \ref{conc_unif_thm} we prove that the number of identifiable admixture components concentrates around its expected value.

In section \ref{choq} we consider inference when the number of identifiable admixture components is equal to $J$, the dimension of the Euclidean space $\mathbb{R}^J$ we are working with, but the admixture components and weights are unknown. In Theorem \ref{thm1}, we use Theorem \ref{choq2} and Glivenko-Cantelli to retrieve the Choquet measure and thus also the identifiable admixture weights.


In section \ref{algorithm}, we use the idea of mixture models based on the extremal set to formulate a novel algorithm that outputs an admixture model composed of only extremal elements, that is, an estimate of the richest cheap model. We state the objective function the algorithm optimizes, and provide a two-stage procedure. We apply this latter to the Associated Press data from the First Text Retrieval Conference (TREC-1), a large collection of terms used in 2246 documents.

Section \ref{concl} concludes our work. In Appendix \ref{app1} we discuss the similarities between de Finetti's theorem and Theorem \ref{choq2}, and we give an approximation of the joint distribution of the admixture components. We also provide the number of extremal elements of the convex hull within a unit simplex having the least amount of vertices.
We prove our results in Appendix \ref{proofs}.

\section{Growth rates for extremal elements and admixture components}\label{mainresults}

\begin{figure}[h!]
\centering
\begin{subfigure}{.5\textwidth}
  \centering
  \includegraphics[width=.7\linewidth]{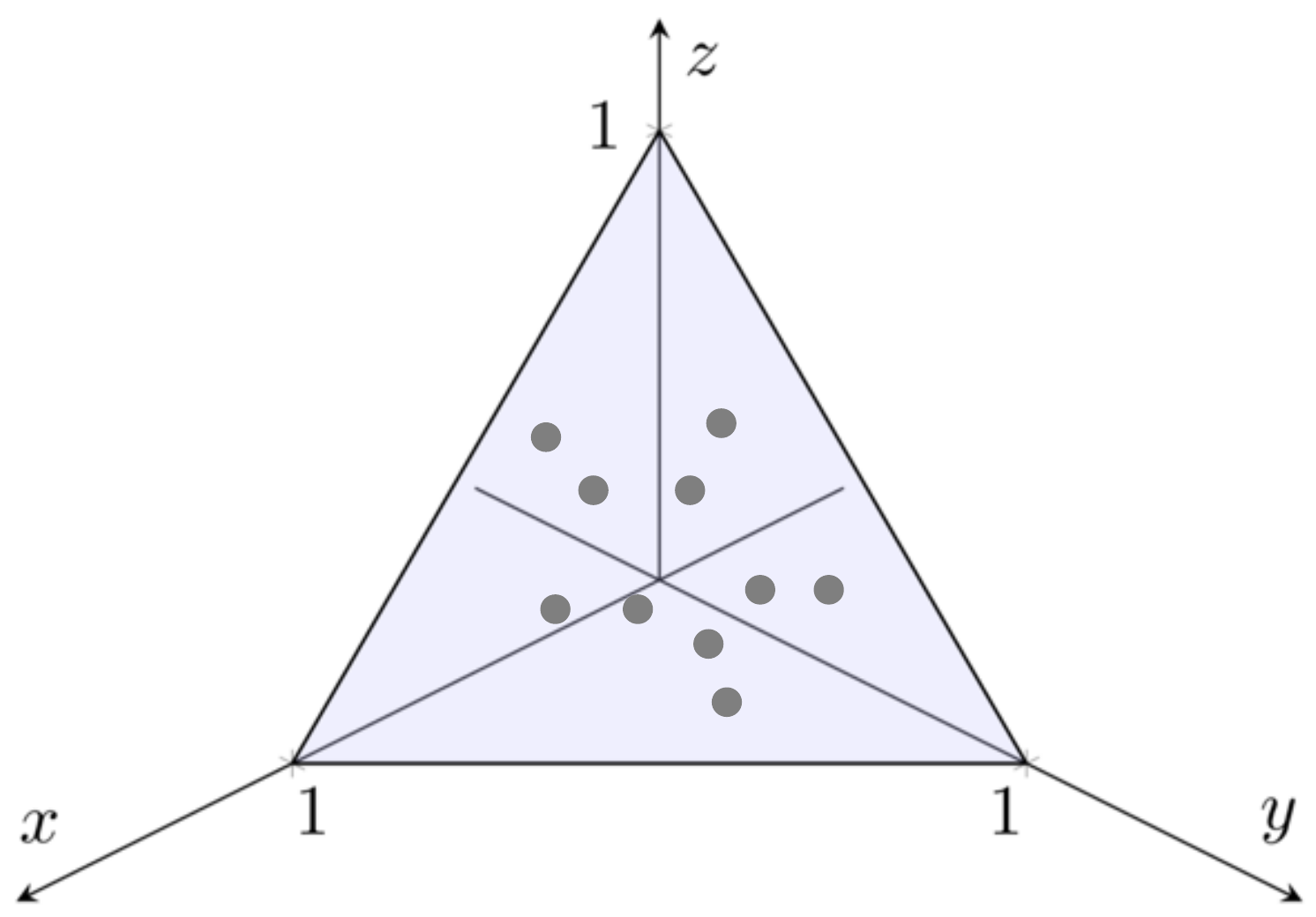}
  \caption{Figure \ref{fig:sub1}}
  \label{fig:sub1}
\end{subfigure}%
\begin{subfigure}{.5\textwidth}
  \centering
  \includegraphics[width=.7\linewidth]{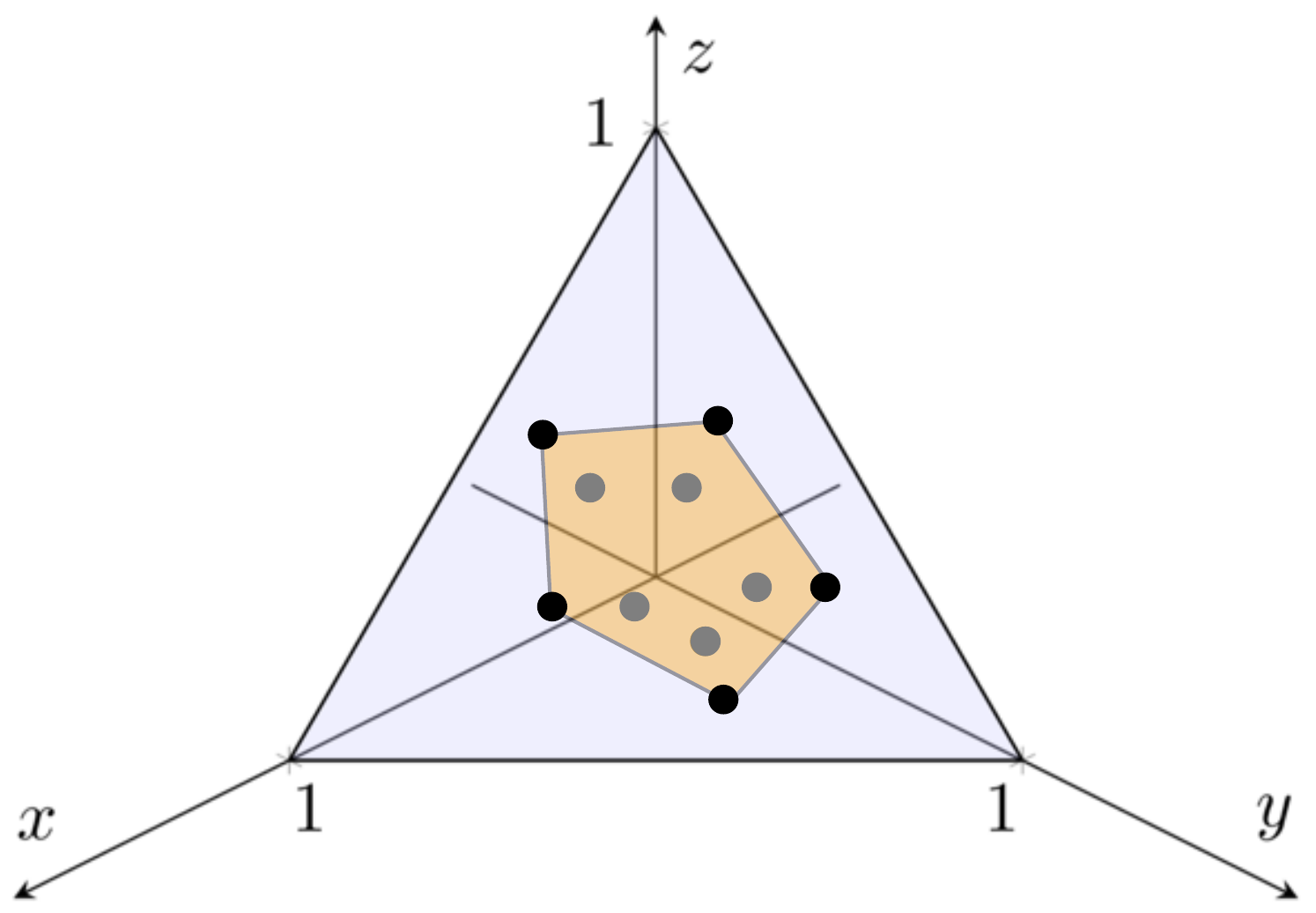}
  \caption{Figure \ref{fig:sub2}}
  \label{fig:sub2}
\end{subfigure}
\caption{Suppose that we are working in $\mathbb{R}^3$. We sample $S_1,\ldots,S_{10}\sim \text{Uniform}(\Delta^2)$ iid; the realizations $s_1,\ldots,s_{10}$ are the ten gray points in the purple unit $2$-simplex in Figure \ref{fig:sub1}. Their convex hull $K_{10}$ is the orange polygon in Figure \ref{fig:sub2}. As we can see, it has five vertices, so $M(10)=5$.}
\label{fig:procedure}
\end{figure}

Suppose the number $M$ of identifiable admixture components is a function $M(n)$ of the amount $n$ of data $x_1,\ldots,x_n$ that we sample. Such function is defined as follows. Let
\begin{equation}\label{eq1_fmm}
S_1,\ldots,S_n \stackrel{iid}{\sim} \mbox{Uniform}(\Delta^{J-1}), \quad J\geq 2,
\end{equation}
and call $K_n:=\text{Conv}(s_1,\ldots,s_n)$, where $s_j$ denotes the realization of $S_j$. In the stochastic convex geometry literature \cite{vu}, $K_n$ is called a \textit{random polytope}. Function $M(\cdot)$ is defined as
$$M:\mathbb{N} \rightarrow \mathbb{N}, \quad n \mapsto M(n):=\#\text{ex}(K_n),$$
that is, $M(n)$ is given by the cardinality of the extremal set of $K_n$. This construction allows us to explicitly relate the field of stochastic convex geometry with applied probability. A simple representation of the procedure to elicit $M(n)$ when $n=10$ and $J=3$ is given in Figure \ref{fig:procedure}.

Before presenting the results in this section, we need to introduce the following geometric concepts.
\begin{itemize}
    \item As pointed out in \cite[Definition 2.1]{ziegler}, in higher-dimensional geometry, the faces of a polytope are features of all dimensions. A face of dimension $i$ is called an \textit{$i$-face}. For example, the polygonal faces of an ordinary polyhedron are $2$-faces. For any $h$-dimensional polytope, we have that $-1 \leq i \leq h$, where $-1$ is the dimension of the empty set. Let us  give a clarifying example. The faces of a cube comprise the cube itself ($3$-face), its facets ($2$-faces), the edges ($1$-faces), its vertices ($0$-faces), and the empty set (having dimension $-1$). Given a generic $h$-dimensional polytope $P$, we denote by $\mathcal{F}_i(P)$ one of its $i$-faces, $i \in \{-1,0,\ldots,h\}$.
    \item We call $\mathscr{F}_i(P)$ the \textit{collection of its $i$-faces}, and ${F}_i(P)$ the \textit{number of its $i$-faces}, that is, ${F}_i(P)=\#\mathscr{F}_i(P)$, for all $i$.
\end{itemize}

In view of the above definition of $0$-faces, we denote by $F_0(K_n)$ the number of extremal elements of $K_n$, so $M(n)=F_0(K_n)$. In the remainder of this section, we keep both notations to highlight the relationship between (stochastic convex) geometry and finite admixture models.

First we find the expected number of identifiable components and we show that it grows at rate $(\log n)^{J-1}$. Equations \eqref{first_theorem} and \eqref{equation_imp_first} are a consequence of \cite[Theorem 6]{reitzner} and \cite[Theorem 5]{barany2}, respectively.

\begin{theorem}{\textbf{(Growth rate of $\mathbb{E}[M(n)]$).}}\label{growth}
Let $K_n :=\text{Conv}({s}_1,\ldots,{s}_{n})$, where ${S}_1,\ldots,{S}_{n}$ are sampled as in \eqref{eq1_fmm}. Then,
\begin{align}\label{first_theorem}
\mathbb{E}[M(n)] = \mathbb{E} \left[F_0(K_n)\right] = \frac{J}{(J+1)^{J-1}}  (\log n)^{J-1} + \mathcal{O}\left((\log n)^{J-2}\log\log n\right),
\end{align}
where $\mathcal{O}$ denotes Bachmann–Landau big-O notation.
In turn, this implies that
\begin{equation}\label{equation_imp_first}
    \lim_{n \rightarrow \infty} (\log n)^{-(J-1)}\mathbb{E} \left[M(n)\right] = \lim_{n \rightarrow \infty} (\log n)^{-(J-1)}\mathbb{E} \left[F_0(K_n)\right] = \frac{J}{(J+1)^{J-1}}  =: c(J).
\end{equation}
\end{theorem}


Then, we see how, for $n$ large enough, the variance $\mathbb{V}[M(n)]=\mathbb{V}[F_0(K_n)]$ of the number of identifiable components can be approximated by $(\log n)^{J-1}$. Equation \eqref{var_eq} is a consequence of \cite[Theorem 1.3]{reitzner_var}.

\begin{theorem}{\textbf{(Approximation of $\mathbb{V}[M(n)]$).}}\label{var_thm}
Let $K_n :=\text{Conv}({s}_1,\ldots,{s}_{n})$, where ${S}_1,\ldots,{S}_{n}$ are sampled as in \eqref{eq1_fmm}. Then, 
\begin{equation}\label{var_eq}
    \mathbb{V}[M(n)]=\mathbb{V}[F_0(K_n)]=\mathcal{O}\left( (\log n)^{J-1} \right).
\end{equation}
\end{theorem}

If $J=3$, we have the following central limit theorem for the number of identifiable admixture components.

\begin{theorem}{\textbf{(CLT for $M(n)$ when $J=3$).}}\label{clt_2}
Let $J=3$ and $K_n :=\text{Conv}({s}_1,\ldots,{s}_{n})$, where ${S}_1,\ldots,{S}_{n}$ are sampled as in \eqref{eq1_fmm}. Then, 
\begin{align}\label{clt_2_eq}
\begin{split}
    \lim_{n\rightarrow\infty} &\sup_{x\in\mathbb{R}} \left| \mathbb{P}\left( \frac{M(n)-\mathbb{E}[M(n)]}{\sqrt{\mathbb{V}[M(n)]}} \leq x \right) - \Phi(x) \right|\\
    &= \lim_{n\rightarrow\infty} \sup_{x\in\mathbb{R}} \left| \mathbb{P}\left( \frac{F_0(K_n)-\mathbb{E}[F_0(K_n)]}{\sqrt{\mathbb{V}[F_0(K_n)]}} \leq x \right) - \Phi(x) \right|=0,
\end{split}
\end{align}
where $\Phi$ denotes the cdf of a Standard Normal distribution.
\end{theorem}

\begin{figure}[h!]
\centering
\begin{subfigure}{.5\textwidth}
  \centering
  \includegraphics[width=.73\linewidth]{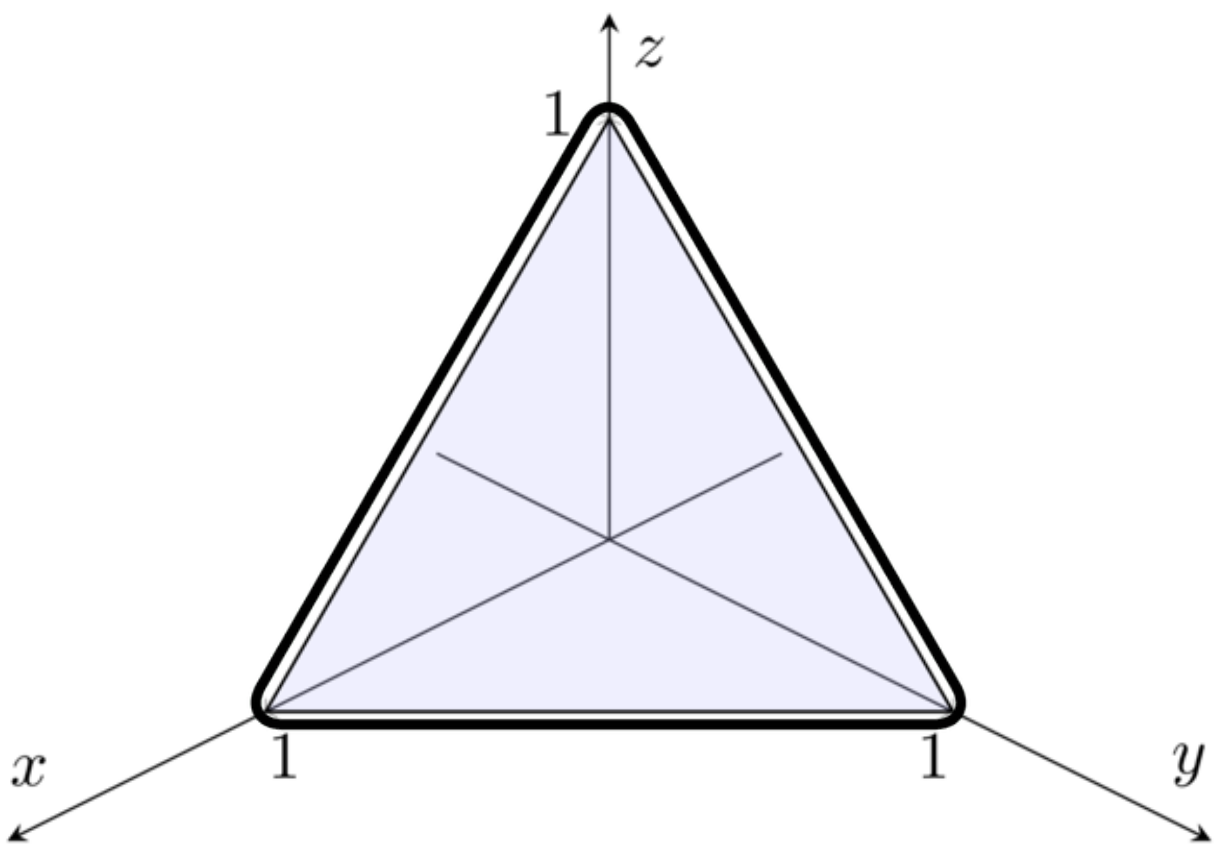}
  \caption{Figure \ref{fig:sub3}}
  \label{fig:sub3}
\end{subfigure}%
\begin{subfigure}{.5\textwidth}
  \centering
  \includegraphics[width=.5\linewidth]{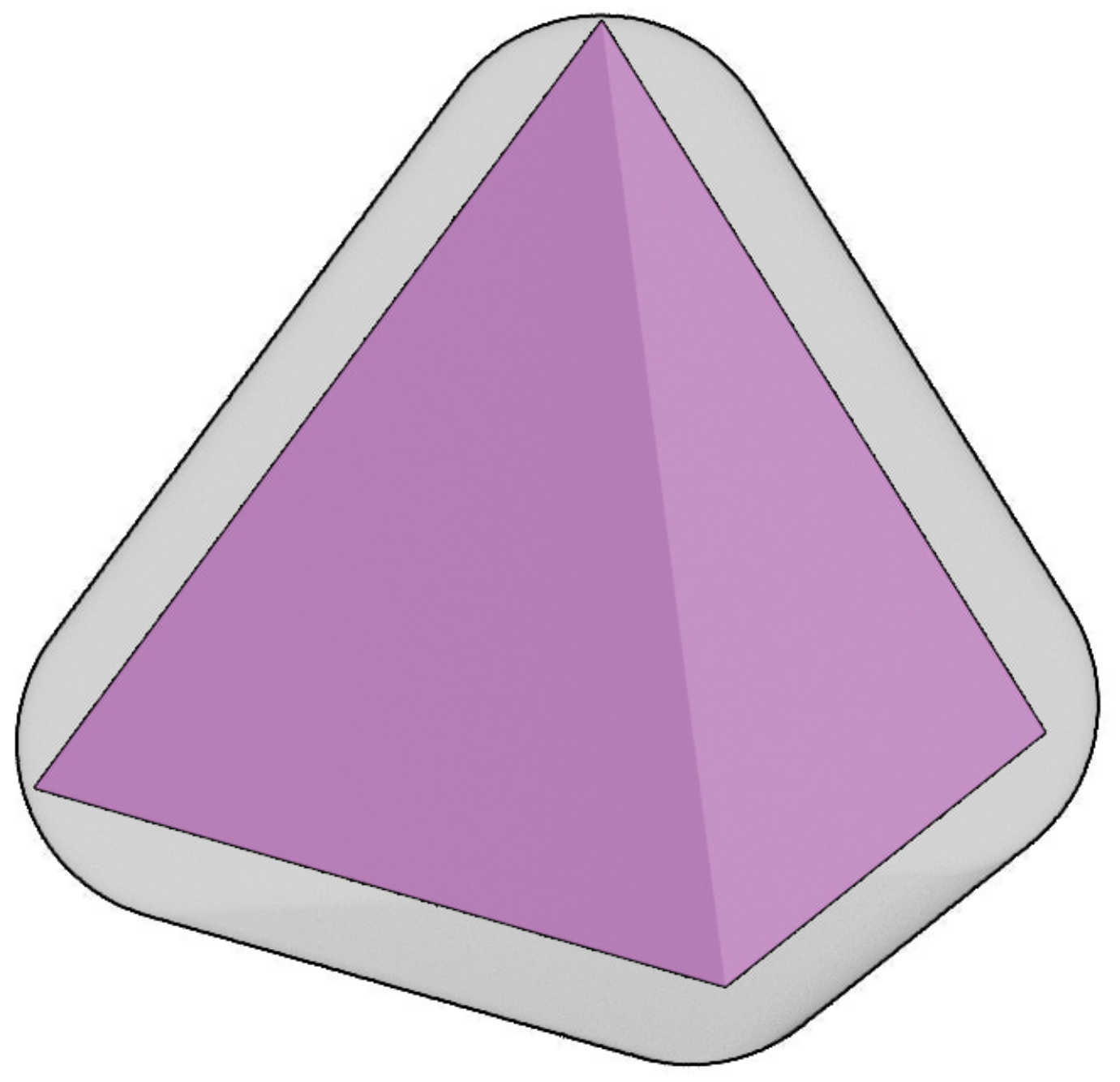}
  \caption{Figure \ref{fig:sub4}}
  \label{fig:sub4}
\end{subfigure}
\caption{For some $\epsilon>0$, $\hat{\Delta}^{2}_\epsilon$ is given by the triangle with round edges in Figure \ref{fig:sub3} containing $\Delta^2$, the unit $2$-simplex in $\mathbb{R}^3$. $\hat{\Delta}^{3}_\epsilon$, instead, is given by the smooth tetrahedron in Figure \ref{fig:sub4} containing $\Delta^3$, the unit $3$-simplex in $\mathbb{R}^4$.}
\label{fig_smooth_simpl}
\end{figure}

Equation \eqref{clt_2_eq} is a consequence of \cite[Corollary 1.2]{pardon}. There, the author conjectures that it holds also for $J> 3$, but to the best of our knowledge such a conjecture has not been proven yet. To overcome this shortcoming, consider the following modification to our setup. Call $\mathbb{K}^2_+(J)$ the space of compact convex sets in $\mathbb{R}^J$, $J\geq 2$, with nonempty interior, boundary of differentiability class $\mathcal{C}^2$, and positive Gaussian curvature. That is, $\mathbb{K}^2_+(J)$ is the space of smooth compact convex sets in $\mathbb{R}^J$; the unit simplex ${\Delta}^{J-1}$ \textit{does not} belong to $\mathbb{K}^2_+(J)$. Pick then any $\epsilon>0$, and call $\hat{\Delta}^{J-1}_\epsilon$ a set in $\mathbb{K}^2_+(J)$ such that ${\Delta}^{J-1} \subset \hat{\Delta}^{J-1}_\epsilon$ and
$$d_H\left( \Delta^{J-1},\hat{\Delta}^{J-1}_\epsilon \right):=\max\left\lbrace{\sup_{x\in \Delta^{J-1}} d_2(x,\hat{\Delta}^{J-1}_\epsilon), \sup_{y\in \hat{\Delta}^{J-1}_\epsilon} d_2({\Delta}^{J-1},y)}\right\rbrace = \epsilon,$$
where $d_H$ denotes the Hausdorff distance and $d_2$ the Euclidean metric. Examples of  $\hat{\Delta}^{2}_\epsilon$ and $\hat{\Delta}^{3}_\epsilon$ are given in Figure \ref{fig_smooth_simpl}. Notice that set $\hat{\Delta}^{J-1}_\epsilon$ always exists, and that it is an $\epsilon$-approximation of ${\Delta}^{J-1}$ belonging to $\mathbb{K}^2_+(J)$. 

Now, sample 
\begin{equation}\label{eq1_hat}
    \hat{S}_1,\ldots,\hat{S}_n \stackrel{iid}{\sim} \text{Uniform}(\hat{\Delta}^{J-1}_\epsilon), \quad J\geq 2,
\end{equation}
and call $\hat{K}_n:=\text{Conv}(\hat{s}_1,\ldots\hat{s}_n)$, where $\hat{s}_j$ is the realization of $\hat{S}_j$, $j\in\{1,\ldots,n\}$. Notice that in \eqref{eq1_hat} we could sample elements that are in $\hat{\Delta}^{J-1}_\epsilon\setminus {\Delta}^{J-1}$, but this happens with a probability that shrinks with $\epsilon$. 
Define then 
$$\hat{M}:\mathbb{N}\rightarrow\mathbb{N}, \quad n \mapsto \hat{M}(n):=\#\text{ex}(\hat{K}_n)=F_0(\hat{K}_n).$$
We can now give a version of Theorem \ref{clt_2} that holds for any $J\geq 2$. Equation \eqref{clt_gen_eq} is a consequence of \cite[Theorems 2, 6]{reitzner1}.

\begin{theorem}{\textbf{(CLT for $\hat{M}(n)$).}}\label{clt_gen}
Let $\hat{K}_n:=\text{Conv}(\hat{s}_1,\ldots\hat{s}_n)$, where $\hat{S}_1,\ldots,\hat{S}_n$ are sampled as in \eqref{eq1_hat}. Then, 
\begin{align}\label{clt_gen_eq}
\begin{split}
    &\left| \mathbb{P}\left( \frac{\hat{M}(n)-\mathbb{E}[\hat{M}(n)])}{\sqrt{\mathbb{V}[\hat{M}(n)]}} \leq x \right) - \Phi(x)\right| \\
    &=\left| \mathbb{P}\left( \frac{F_0(\hat{K}_n)-\mathbb{E}[F_0(\hat{K}_n)])}{\sqrt{\mathbb{V}[F_0(\hat{K}_n)]}} \leq x \right) - \Phi(x)\right| = \mathcal{O}\left( n^{-\frac{1}{2(J+1)}}\left(\log n\right)^{2+\frac{2}{J+1}} \right).
\end{split}
\end{align}
\end{theorem}

Then we prove that the number of identifiable admixture components concentrates around its expected value. Equation \eqref{conc_unif_eq} is a consequence of \cite[Theorem 2.11, Section 7]{vu}.

\begin{theorem}{\textbf{(Concentration inequality for $\hat{M}(n)$).}}\label{conc_unif_thm}
    Let $\hat{K}_n:=\text{Conv}(\hat{s}_1,\ldots\hat{s}_n)$, where $\hat{S}_1,\ldots,\hat{S}_n$ are sampled as in \eqref{eq1_hat}. Then, there are fixed positive constants $c,\Xi,\epsilon_0$ such that for any $\epsilon\in(0,\epsilon_0]$, $V\geq \Xi n^{\frac{J-1}{J+1}}$, $C\geq n\epsilon$, and $\lambda\in(0,\frac{V}{4C^2})$, the following holds
    \begin{align}\label{conc_unif_eq}
        \begin{split}
            \mathbb{P}&\left( \left| \hat{M}(n) -\mathbb{E}[\hat{M}(n)] \right| \geq \sqrt{\lambda V} \right)\\
            &=\mathbb{P}\left( \left| F_0(\hat{K}_n) -\mathbb{E}[F_0(\hat{K}_n)] \right| \geq \sqrt{\lambda V} \right) \leq 2 \exp(-\lambda/4)+\exp(-c\epsilon n) + \exp\left(-cn^\frac{J-1}{3J+5}\right).
        \end{split}  
    \end{align}
\end{theorem}

\subsection{From the Uniform to the general case}\label{from-to}
In this section we generalize Theorem \ref{growth} to the case where $S_1,\ldots,S_n$ are iid samples from a generic distribution $G$.\footnote{That is, $G$ is any distribution that is absolutely continuous with respect to the Lebesgue volume measure of the simplex.} We defer the study of the non-iid case to future work. We ask ourselves whether the asymptotic growth function of the expected number of identifiable components based on draws from the uniform distribution can inform us about the growth rate based on draws from a generic distribution. 

Suppose $S_1,\ldots,S_n$ are now sampled iid from a generic distribution $G$ on $\Delta^{J-1}$. Call then $\breve{K}_n:=\text{Conv}(s_1,\ldots,s_n)$ and define
$$T:\mathbb{N} \rightarrow \mathbb{N}, \quad n \mapsto T(n):=\#\text{ex}(\breve{K}_n)=F_0(\breve{K}_n),$$
so we assume that the number of identifiable admixture component corresponds to the number of vertices of $\breve{K}_n$. In a way, we can see $T(n)$ as a ``generalization'' of $M(n)$; of course, $T(n)$ and $M(n)$ may be different.





\begin{remark}
In this work, we only consider the case of $M(n),T(n)\geq J$ (of course, $J \geq 2$). This because, if that were not the case, our convex hull would be a (proper) subset of a smaller-dimensional Euclidean space than $\mathbb{R}^J$. 
It is worth to notice, though, that {\em if $n \geq J$ and is sufficiently large, then $M(n),T(n)\geq J$ holds almost surely}. Let us be more precise. In Theorem \ref{growth}, we show that $\mathbb{E}[M(n)]$ grows as $(\log n)^{J-1}$, as $n\rightarrow \infty$. By concentration of measure arguments, then, it follows that the probability of having only \( J \) extreme points is vanishingly small. 
More formally, if we assume \( n \geq J \), it follows that for sufficiently large \( n \), the expected number of extreme points satisfies
\[
\mathbb{E}[M(n)] \approx c(J) (\log n)^{J-1} \geq J.
\]
Thus, by Markov's inequality and standard deviation bounds, we have that

\[
\mathbb{P}(M(n) < J) \rightarrow 0 \quad \text{as } n \rightarrow \infty.
\]
That is,

\[
\mathbb{P}(M(n) \geq J) \rightarrow 1 \quad \text{as } n \rightarrow \infty.
\]
Therefore, almost surely, as we sample more and more points uniformly from $\Delta^{J-1}$, the number $M(n)$ of extreme elements of $K_n$ is greater than \( J \). A similar argument holds for $T(n)$ as a consequence of Theorem \ref{extrema2}. 

This argument shows that, while it is necessary to assume that $M(n),T(n) \geq J$ if we consider a fixed generic sample size $n$, this condition is satisfied almost surely when $n \rightarrow \infty$.
\end{remark}

The next result, Theorem \ref{extrema2}, states the following. Up until the $(N-1)$-th data point, the expected number of identifiable components $\mathbb{E}[T(n)]$ in the more general case can take on any possible real value. From the $N$-th observation onward, though, it must be in a fixed (possibly highly nonlinear) relationship $\gamma_n$ with $\mathbb{E}[M(n)]$. If this happens, we are able to relate their growth rates. Such an assumption is made primarily for mathematical convenience: with it, deriving the result in equation \eqref{gen_gr_rate} becomes relatively easy. In the future, we plan to relax this assumption in order to derive a more general result.

\begin{theorem}{\textbf{(Growth rate of $\mathbb{E}[T(n)]$).}}\label{extrema2}
Call $(\gamma_n)$ a sequence in $\mathbb{R}^\mathbb{N}$ for which $0$ is not an accumulation point, and let $\mathbb{E}[T(n)]=g_n(\mathbb{E}[M(n)])$, where $g_n$ is a functional on $\mathbb{R}$ that depends on $n$. Then, if there exists $N\in\mathbb{N}$ such that for all $n\geq N$, $g_n(\mathbb{E}[M(n)])=\gamma_n \mathbb{E}[M(n)]$, we have that
\begin{equation}\label{gen_gr_rate}
\lim_{n \rightarrow \infty} \frac{1}{{\gamma_n} (\log n)^{J-1}} \mathbb{E}[T(n)] =  c(J).
\end{equation}
\end{theorem}

In \eqref{gen_gr_rate}, $c(J)$ is a quantity that only depends on the dimension $J$ of the Euclidean space $\mathbb{R}^J$ that we consider, and it is the same as in \eqref{equation_imp_first}. The following corollary is a direct consequence of Theorem \ref{extrema2}.

\begin{corollary}\label{extrema2_cor}
Suppose the assumptions of Theorem \ref{extrema2} hold. Then, if there exists a sequence $\varpi_n \in \mathbb{R}^\mathbb{N}$ such that $\gamma_n=\mathcal{O}(\varpi_n)$, then the growth rate of $\mathbb{E}[T(n)]$ is $\varpi_n(\log n)^{J-1}$.
\end{corollary}

\begin{remark}
It is immediate to see that there is a universal upper bound for the Euclidean distance between two points in a unit simplex: for all $x,y \in \Delta^{J-1}$, $d_2(x,y) \equiv \|x-y\| \leq 2$, where $\|\cdot\|$ denotes the Euclidean norm. This gives us an interesting result: the Hausdorff distance between $K_n$ and $\breve{K}_n$ has a universal upper bound as well. Indeed,  
\begin{align*}
    d_H({K}_n,\breve{{K}}_n)&=  \max\left\lbrace{ \sup_{{x} \in {K}_n} d_2 (x,\breve{{K}}_n), \sup_{y \in \breve{{K}}_n} d_2({K}_n,y)}\right\rbrace  \\
    &=\max\left\lbrace{ \sup_{{x} \in {K}_n} \inf_{y \in \breve{{K}}_n}d_2 (x,y), \sup_{y \in \breve{{K}}_n} \inf_{x \in {K}_n} d_2(x,y)}\right\rbrace \leq 2.
\end{align*}

Notice also that if instead of requiring $\mathbb{E}[T(n)]=\gamma_n \mathbb{E}[M(n)]$, for all $n\geq N$,   we are willing to make the slightly stronger assumption that for all $n\geq N$, $T(n)=\rho_n  M(n)$, $(\rho_n) \in \mathbb{R}^\mathbb{N}$ possibly different from $(\gamma_n)$, then we retrieve Theorem \ref{clt_2} for $T(n)$.\footnote{Of course we need to assume that for $(\rho_n)$, $0$ is not an accumulation point.} This because, since $F_0(K_n) = M(n)$, we have that 
$$\frac{T(n) -\mathbb{E}[T(n)]}{\sqrt{\mathbb{V}[T(n)]}} = \frac{\rho_n F_0(K_n) - \mathbb{E}[\rho_n F_0(K_n)]}{\sqrt{\mathbb{V}[\rho_n F_0(K_n)]}} = \frac{F_0(K_n) -\mathbb{E}[F_0(K_n)]}{\sqrt{\mathbb{V}[F_0(K_n)]}},$$
and so Theorem \ref{clt_2} follows. A similar argument allows us to retrieve Theorem \ref{clt_gen} when we work with $\hat{K}_n$ instead of ${K}_n$.
\end{remark}



\section{Choquet measure and identifiable admixture weights}\label{choq}

In this section we build a bridge between finite admixture models and Choquet theory. We show how, thanks to a uniqueness result by Gustave Choquet (i.e. Theorem \ref{choq2}), the classical Glivenko-Cantelli's theorem can be used to retrieve the identifiable admixture weights. 


By Theorem \ref{choq2}, we have that for every element $p$ in a simplex $C$, there exists a unique measure -- that we call the Choquet measure associated with $p$, and denote by $\nu_p$ -- supported on the extremal elements $E=\text{ex}(C)$ such that $p=\sum_{e \in E} e \cdot \nu_p(e)$.\footnote{We write $\nu_p(e)$ in place of $\nu_p(\{e\})$ for notational convenience. We stick to this abuse of notation for the rest of the paper.} In our analysis, $p$ corresponds to $\pi$, the elements $e$ in $E=\text{ex}(C)$ correspond to the identifiable $f_\ell$'s, and the $\nu_p(e)$'s correspond to the weights of the identifiable $f_\ell$'s, that is, $\nu_{\pi}(f_\ell)=\varphi_{\ell}$, for every identifiable $f_\ell$.

In Theorem \ref{thm1}, we show that if we only assume that the number $M$ of components is known and equal to $J$ -- the dimension of the Euclidean space $\mathbb{R}^J$ we are working with -- we can use Glivenko-Cantelli to retrieve $\nu_{\pi}$. The $\varphi_{\ell}$'s represent the weights of the richest cheap finite admixture model representing $\pi$, 
so $\pi=\sum_{\ell=1}^M f_\ell \nu_{\pi}(f_\ell)$. 


\begin{remark}
Recall that $\pi=\sum_{\ell=1}^M \varphi_{\ell} f_\ell =\sum_{\ell=1}^L \phi_{\ell} f_\ell $, where we labeled the unidentifiable components as $f_{M+1},\ldots,f_L$, $M \leq L$. This is without loss of generality.
\end{remark}

\begin{figure}[h!]
\centering
\includegraphics[width=.45\textwidth]{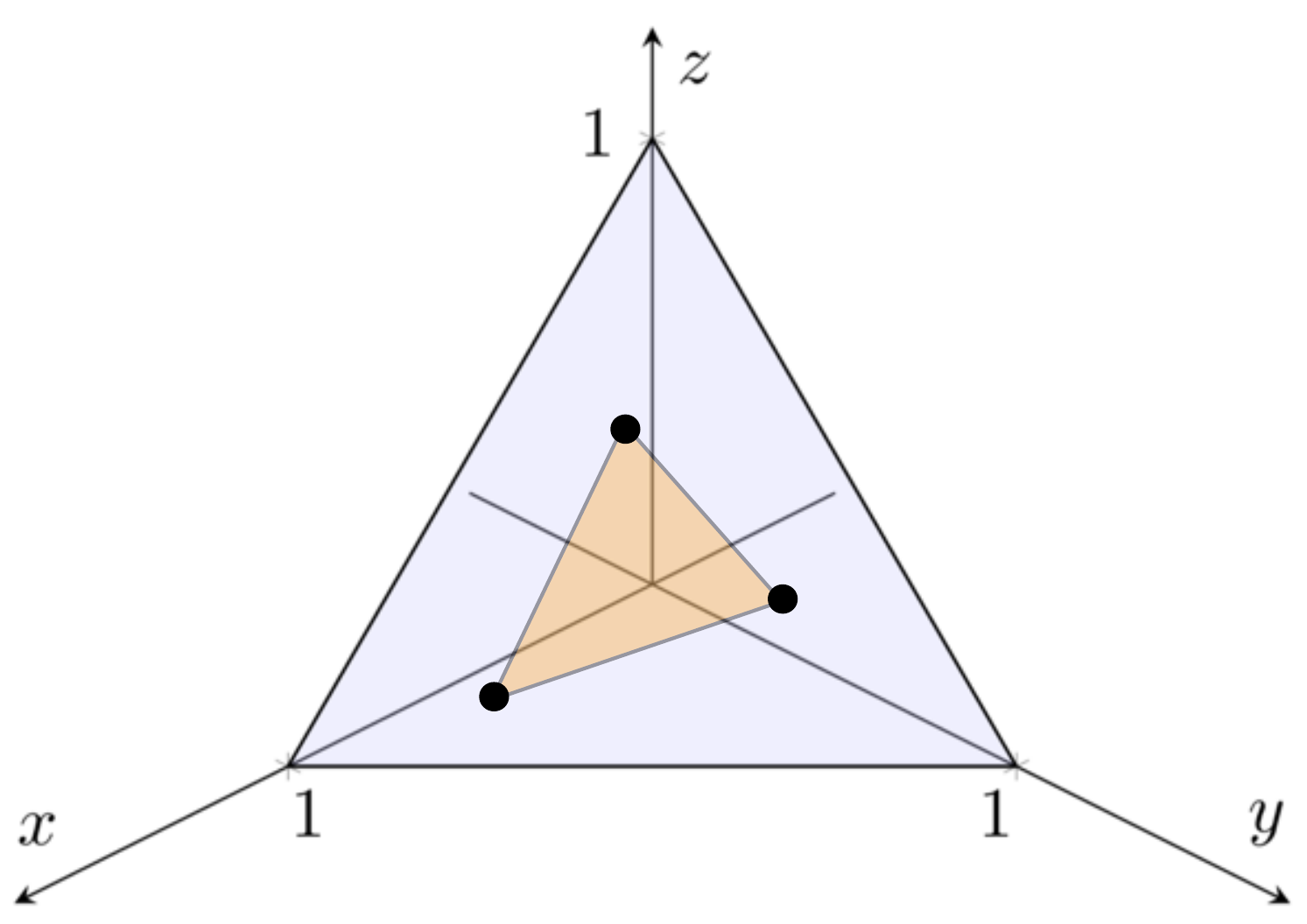}
\caption{A triangular-shaped convex hull within the unit $2$-simplex in $\mathbb{R}^3$. It is a simplex because the number of its vertices coincides with the dimension of the Euclidean space.}
\label{fig2}
\centering
\end{figure}

Let us denote by 
$$\mathcal{K}_M:=\text{Conv}(f_1,\ldots,f_M)=\text{Conv}(f_1,\ldots,f_L)$$ 
the convex hull generated by the $M$ identifiable admixture components, and assume $M=J$, so that $\mathcal{K}_M$ is a {\em simplex}. Recall that the latter is defined as $\mathcal{K}_M \coloneqq \{\pi=\sum_{\ell=1}^M \xi_\ell f_\ell : \sum_{\ell=1}^M \xi_\ell =1 \text{, } \xi_\ell \geq 0 \text{ } \forall \ell \in \{1,\ldots,M\} \text{, and } M=J\} \subset \mathbb{R}^J$, while the {\em $(J-1)$-unit simplex} satisfies the extra condition that the $f_\ell$'s are the $J$ standard unit vectors in $\mathbb{R}^J$. An example of a simplex within the unit $2$-simplex in $\mathbb{R}^3$ is given in Figure \ref{fig2}.  Our first goal is to learn about distributions supported on the extremal elements of $\mathcal{K}_M$, $E_M:=\text{ex}(\mathcal{K}_M)=\{f_1,\ldots,f_M\}$.

Since $\Delta^{J-1}$ is locally convex, and $\mathcal{K}_M \subset \Delta^{J-1}$ is a metrizable compact convex set, then thanks to Theorem \ref{choq2}, we know that for every $\pi \in \mathcal{K}_M$, there exists a unique probability measure $\nu_{\pi}$ supported on $E_M$ such that $\pi=\sum_{f_\ell \in E_M} f_\ell \cdot \nu_{\pi}(f_\ell)$. We formalize this statement in the next proposition.



\begin{proposition}{\textbf{(Choquet uniqueness for admixture models)}.}\label{simplex_conseq}
If $\mathcal{K}_M$ is a simplex, then every element $\pi \in \mathcal{K}_M$ can be represented by a unique measure $\nu_{\pi}$ (the Choquet measure representing $\pi$) supported on $E_M$.
\end{proposition}
For every identifiable admixture component $f_\ell\in E_M$, the Choquet measure $\nu_{\pi}$ gives the corresponding admixture weight, that is, $\nu_{\pi}(f_\ell)=\varphi_{\ell}$.

We are now ready for the main result of this section. Suppose that $\mathcal{K}_M$ is a simplex, so that $M=J$, and call $E_M$ the set of extremal elements of $\mathcal{K}_M$. We denote by $e_j$ a generic element of $E_M$, that is, a generic identifiable admixture component. In mathematical terms, we write $e_j\in\{f_1,\ldots,f_M\}=:E_M$, for all $j\in\mathbb{N}$. 

Consider an iid sample from the Choquet measure $\nu_{\pi}$, that is, $e_1,\ldots,e_K \sim \nu_{\pi}$ iid. Because every $e_k$, $k\in \{1,\ldots,K\}$, corresponds to an identifiable admixture element, we can write $e_k=f_{\ell_k}$, where label $\ell_k$ belongs to $\{1,\ldots,M\}$, for all $k\in\{1,\ldots,K\}$. Now, call $\zeta$ the distribution of the labels; we immediately notice that there is a one-to-one correspondence between $\zeta$ and $\nu_{\pi}$ since $\zeta(\ell)=\nu_{\pi}(f_{\ell})$, for all $\ell\in\{1,\ldots,M\}$. Call then $F_\zeta$ the cdf of $\zeta$, and denote by $\ell_1,\ldots,\ell_K\sim \zeta$ the iid sample of labels from $\zeta$ that corresponds to the iid sample $e_1,\ldots,e_K$ from the Choquet measure $\nu_{\pi}$. For all $x\in\mathbb{R}$, denote by
$$F_K(x):=\frac{1}{K}\sum_{j=1}^K \mathbb{I}_{[\ell_j,\infty)}(x)$$
the empirical cdf of $\ell_1,\ldots,\ell_K$, where $\mathbb{I}_A(x)$ stands for the indicator function of $x$ belonging to a generic set $A$. Then, we have the following. Equation \eqref{gk_cvg} is a consequence of the Glivenko-Cantelli's Theorem \cite{vdv_asy}. 

\begin{theorem}{\textbf{($F_K$ retrieves $F_\zeta$)}.}\label{thm1}
    Let $\mathcal{K}_M$ be a simplex in $\Delta^{J-1}$, call $E_M:=\text{ex}(\mathcal{K}_M)$ the set of its extremal elements, so that $M=J$. Suppose that, in the notation introduced above, $\ell_1,\ldots,\ell_K\sim \zeta$ iid. Then, we have that 
\begin{equation}\label{gk_cvg}
    \sup_{x\in\mathbb{R}} \left| F_K(x)-F_\zeta(x) \right| \xrightarrow[K\rightarrow\infty]{a.s.} 0.
\end{equation}
In addition, pick any $0 < \eta < 0.5$. Then, $\sup_{x\in\mathbb{R}} \left| F_K(x)-F_\zeta(x) \right| = o \left( \frac{1}{K^\eta} \right)$ a.s.
\end{theorem}

The idea in Theorem \ref{thm1} is that if we keep observing iid samples $e_k$'s from the Choquet measure $\nu_{\pi}$ -- that is, if we are able to observe identifiable admixture components $\{f_1,\ldots,f_M\}$ sampled according to their true weights $\{\nu_{\pi}(f_1),\ldots,\nu_{\pi}(f_M)\}\equiv\{\varphi_{1},\ldots,\varphi_{M}\}$ -- then their empirical cdf recovers the cdf of $\nu_{\pi}$. In turn, this gives us the identifiable admixture weights, since -- as we saw before -- $\nu_{\pi}(f_\ell)=\varphi_{\ell}$, for all $\ell\in\{1,\ldots,M\}$.

\section{A procedure to estimate the richest cheap model}\label{algorithm}


In this section, we present a sparse-factor-analysis-inspired algorithm that estimates the parameter of the richest cheap model introduced in section \ref{intro}, and we apply it to the well-know TREC-1 document-term matrix dataset \cite{harman}. 

Consider an admixture model $X_i \sim \text{Mult}(\pi_i)$, $\pi_i = \sum_{\ell=1}^L \phi_{i,\ell} f_\ell$, $i\in\{1,\ldots,n\}$, where the $X_i$'s are independent. There are two sets of parameters: 
\begin{enumerate}
\item the mixing weights for each individual, that can be arranged in an $n\times L$ matrix $\Phi$ whose entries $\phi_{i,\ell}$ represent the probability that the
$i$-th sample is drawn from the $\ell$-th component. In turn, each row of $\Phi$ is the mixture vector of the $i$-th observation $\phi_i = (\phi_{i,1},\ldots,\phi_{i,L})$;
\item the probability vectors parameterizing each admixture component, which we can write as an $L\times J$ matrix $F$ whose $\ell$-th row is $f_\ell$.
\end{enumerate}
The relation between admixture modeling and sparse factor analysis has been  
explored in detail in \cite{engel}. There, conditions are provided when sparse factor analysis and LDA have very similar results, and the implications for population genetics are discussed. The key insight in \cite{engel} is that given an $n\times J$ observation matrix $X$ (whose $i$-th row is $X_i$) from a multinomial admixture model, learning an admixture model amounts to the following minimization procedure
\begin{equation}\label{min_probl}
\min_{F,\Phi}\| \mathbb{E}[X] -  \Phi F\|^2.
\end{equation}
The sparse factor analysis framework can be summarized as minimizing \eqref{min_probl} with the constraint that many of the elements of $\Phi$ will be zero, or that every observation is a sparse combination of each component. The spirit behind the algorithm proposed in this section is to think of sparsity as the extremal set: we want to find a set of components that are extremal yet still accurately solves the above minimization. We first state the likelihood for the admixture model, assuming a maximum of $L$ components,
$$\mathcal{L}(X_1,...,X_n; \{\phi_1,...,\phi_n\}, \{f_1,...,f_L\}) = \prod_{i=1}^n \mbox{Mult}\left(\pi_i =  \sum_{\ell=1}^L \phi_{i,\ell} f_\ell  \right).$$
A maximum likelihood estimator (MLE) for for the above model is
\begin{equation}
\label{mixtureml}
 \{\{\hat{\phi}_1,...,\hat{\phi}_n\}, \{\hat{f}_1,...,\hat{f}_L\}\} \in \argmax_{\{\phi_1,...,\phi_n\}, \{f_1,...,f_L\}}  \mathcal{L}(X_1,...,X_n; \{\phi_1,...,\phi_n\}, \{f_1,...,f_L\}).
 \end{equation}
A notion of sparsity related to the sparse factor analysis framework is to maximize the likelihood subject to the constraint that components are identifiable, that is, no component can be represented as a convex combination of other components. We consider a procedure that maximizes the following objective function
\begin{equation} \label{opt_probl}
\begin{aligned}
& \underset{ \mathcal{I},\{\phi_1,...,\phi_n\}, \{f_k\}_{k \in \mathcal{I}}} {\text{argmax} }
& &   \prod_{i=1}^n \mbox{Mult}\left(\pi_i =  \sum_{k \in \mathcal{I}} \phi_{i,k} f_k  \right). \\
& \text{subject to}
& & f_k \not \in \mbox{Conv}(f_{\mathcal{I} \setminus \{k\}}) , \; \forall k \in \mathcal{I} ,
\end{aligned}
\end{equation}
where $\mathcal{I}$ is a subset of the set $\{1,...,L\}$, and it is the collection of the indices of the extremal set. Constraint  $f_k \not \in \mbox{Conv}(f_{\mathcal{I} \setminus \{k\}})$, for all $k \in \mathcal{I}$, ensures that no admixture component is contained in the convex combination of the others. Notice that the cardinality of $\mathcal{I}$ represents the number of components $M$ of the richest cheap model in \eqref{rcm}.

The maximization specified by equation \eqref{opt_probl} is non-convex and finding the global optima is difficult; we propose a two-step procedure to solve it.

\begin{algorithm}
\caption{Approximating the parameters of the richest cheap model}\label{alg:2step}
\begin{algorithmic}
\State \textbf{Step 0} Initialize $t=0$ and set the initial number of components $L_0 \in\mathbb{N}$
\Do:
\State \textbf{Step 1} Compute an MLE as in \eqref{mixtureml} so to obtain the estimated parameters $\{\{\hat{\phi}_1,...,\hat{\phi}_n\}, \{\hat{f}_1,...,\hat{f}_{L_{t}}\}\}$

\State \textbf{Step 2} $L_{t+1}:=\#\text{ex}(\text{Conv}(\{\hat{f}_1,...,\hat{f}_{L_{t}}\})$
\doWhile{$L_{t+1} < L_{t}$} \Comment{\% Call $L_T$ the number of components that exits the loop \%}
\State \Return parameters $\{\{\hat{\phi}_1,...,\hat{\phi}_n\}, \{\hat{f}_1,...,\hat{f}_{L_T}\}\}$
\end{algorithmic}
\end{algorithm}

The parameters returned by Algorithm \ref{alg:2step} are estimates of the parameters of the richest cheap model. Notice that computing the convex hull is evocative of the Choquet procedure described in section \ref{choq}. If $L_T=J$ (where $J$ is the dimension of the Euclidean space we work with), we have an estimate $\hat{\nu}_{\pi_i}$ of the Choquet measure for $\pi_i$, since the geometric representation of the richest cheap model is a simplex. 
We have the following important result.

\begin{theorem}\label{algo_cons}
    Call $M^\star$ the true number of components of the richest cheap finite admixture model. Then, if $L_0 \geq M^\star$ and the MLE in \eqref{mixtureml} is consistent,\footnote{Conditions for the MLE to be consistent can be found in \cite{wasserman_2}.} we have that $L_T$ is a consistent estimator for $M^\star$.
\end{theorem}

We leave for future work to verify whether the following plausible conjecture holds. Suppose that our model is ``sparse enough'' as in \eqref{opt_probl}. That is, suppose we constrain the admixture elements to be identifiable. Then, if the MLE in \eqref{mixtureml} is consistent, it should be verified that the collection $\{\hat{f}_1,...,\hat{f}_{L_T}\}$ of identifiable components estimated by Algorithm \ref{alg:2step} is a consistent estimator for the collection  $\{{f}^\star_1,...,{f}^\star_{M^\star}\}$ of true identifiable admixture components.

\begin{remark}\label{rem_l0}
    One thing left to discuss before applying our algorithm to a dataset is how to choose $L_0$. The objective function \eqref{mixtureml} is highly non-concave and, in general, optimization algorithms
for finding the optimal solutions to \eqref{mixtureml} can have sub-linear convergence
rates. This is due to the fact that function \eqref{mixtureml} is locally weakly concave around
the optimal solutions (see \cite{dwivedi} for a more detailed discussion). Hence, $L_0$ cannot just be chosen to be arbitrarily large, as the performance of the optimization algorithms can be strongly
affected. At the same time, it cannot be too small, otherwise our algorithm would not be able to capture the underlying complexity associated with the data at hand, and we would also risk not to meet the $L_0 \geq M^\star$ condition of Theorem \ref{algo_cons}. We select $L_0$ via an ``educated guess'' coming from the exploratory data analysis part of our study or from previous results on the same or similar datasets. In the future we will study a more formal way of coming up with a value for $L_0$. 

\end{remark}


We applied our two-step procedure to a well studied dataset \cite{bail} which is a document-term matrix consisting of term frequencies of $10473$ terms in $2246$ documents collected from Associated Press documents \cite{harman}. We used the latent Dirichlet allocation (LDA) function in the R package  \texttt{topicmodels} \cite{topicmodels} to compute the MLE and the convex hull function in the R  package \texttt{geometry} to compute the convex hull. 

The number of topics obtained in previous studies on the same dataset is between $9$ and $12$ \cite{bail,hou}. For this reason, we ran our algorithm three times on the document-term matrix setting $L_0$ equal to $12$, $25$, and $50$. $L_0=12$ seems the proper educated guess, while $L_0=25$ and $L_0=50$ are safety checks; starting with a higher value of $L_0$ can cause our algorithm to incur problems, as discussed in Remark \ref{rem_l0}.


Computing the convex hull over the full topic frequency vectors -- elements belonging to simplex $\Delta^{2245}$ -- is prohibitive and also does not make sense when the number of topics are less than $2245$. We used principal components analysis (PCA) to project the frequency vectors of the topics onto a lower dimensional space and then computed the convex hull of the projections. We used a simple scree plot 
to notice that $3-5$ dimensions are sufficient to capture about $40\%$ of the variation when we carry out our analysis specifying $50$ initial topics. If we choose $L_0<50$, we have that $3-5$ dimension explain more than $40\%$ of the variation. We only need to compute the number of extremal elements of the convex hull and not the extremal elements themselves in our procedure, so it suffices to compute the convex hull in the low dimensional space.

\begin{figure}[h!]
\centering
\includegraphics[width=\textwidth]{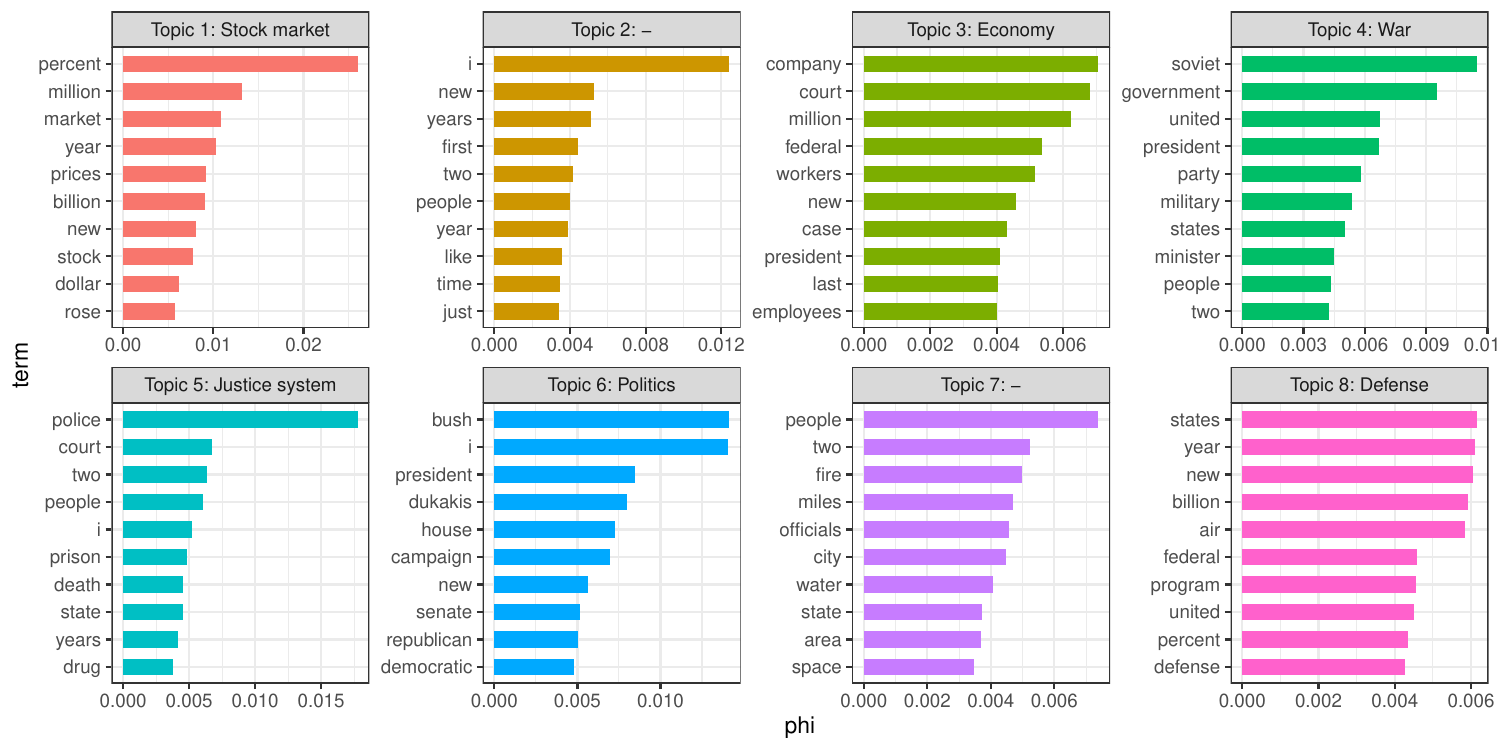}
\caption{Parameters estimated by our algorithm. For every topic, we show the $10$ terms having the highest probability of being generated from that topic. As we can see, Topics 2 and 7 are hard to interpret, a common issue with topic models \cite{fung}.}
\label{fig_algo}
\centering
\end{figure}

Given the results in the PCA step, we projected down to $5$ dimensions. The number of extremal elements -- i.e. the number of topics -- we obtained were  $8$, $8$, and $9$ having initialized $L_0$ to $12$, $25$, and $50$, respectively. The number of topics of the richest cheap model seems to be $8$, that is, a mixture of $8$ multinomials appears to be the model that captures the complexity in our dataset using the smallest number of components. The estimated topics $\ell\in\{1,\ldots,8\}$, together with the $10$ terms $i$ having the highest estimated probability $\hat{\varphi}_{i,\ell}$ of being generated from topic $\ell$, are reported in Figure \ref{fig_algo}. We do not report the estimated topic frequency vectors $\hat{f}_\ell$, $\ell\in\{1,\ldots,8\}$, because of their high dimension, $J=2246$. Recall that for all $\ell\in\{1,\ldots,8\}$, we can write $\hat{f}_\ell$ as $\hat{f}_\ell=(\hat{f}_{\ell,1},\ldots,\hat{f}_{\ell,J})^\top$, where $\hat{f}_{\ell,j}$ represents the estimated probability of topic $\ell$ being featured in document $j$.  

Let us add a brief discussion. In this section we provide a simple application of Algorithm \ref{alg:2step} to show that the latter is actually implementable. A more in-depth discussion on how it performs in comparison to other methods is postponed to future research. We conjecture that Algorithm \ref{alg:2step} finding a smaller number of topics, compared e.g. to \cite{bail,hou}, is due to the fact that these latter look for the ``most interpretable'' solution, that is, the one whose topics are better interpretable by a human user. On the contrary, we explicitly look for the model involving the smallest possible number of admixture components, which is able to fully capture the complexity of the data at hand. That is, we look for the richest cheap model. Another reason why we find a smaller number of topics may linked to the fact that we adopt PCA. Once again, all of these questions are left for future research.

\section{Conclusion}\label{concl}
In this paper there are two key ideas. The first one is that we can use techniques from stochastic convex geometry on the growth rate of the expected number of extremal elements of random polytopes to provide insights into the asymptotic growth rate of the expected number of identifiable admixture components. We prove that the expected number of identifiable components grows at rate $(\log n)^{J-1}$ where $J$ is the dimension of the Euclidean space we work with.
We also show that the number of identifiable admixture components concentrates around its expected value, and we provide a central limit theorem for its distribution. 
The other key concept is that we can retrieve the identifiable admixture weights using techniques from Choquet theory. In particular, we show that if the convex hull $\mathcal{K}_M$ generated by the identifiable admixture components is a simplex, we can apply Glivenko-Cantelli to recover the identifiable admixture weights. We also give an algorithm to estimate the richest cheap finite admixture model. 


An interesting open question is whether there are other instances in (Bayesian) inference where coupling results from stochastic (convex) geometry with results from Choquet theory allows to develop novel analyses, insights, models, or algorithms. For example, studying the properties of an apeirogon -- a polytope with infinitely many sides -- could give us insights on infinite mixture models. Another research direction for the future regards a formal procedure to initialize $L_0$ in Algorithm \ref{alg:2step}. A promising way to tackle this issue is to use the backward induction approach of \cite[Section 1.2]{peskir}, where the authors find an optimal stopping time for a given expected utility maximization problem. In our framework, such optimal stopping time could be interpreted as number $L_0\in\mathbb{N}$ that strikes a balance between being not too large, so to avoid the problems highlighted by \cite{dwivedi}, and not too small, so to satisfy the assumption of Theorem \ref{algo_cons}. We also plan to further generalize the results in section \ref{from-to}.

\section*{Acknowledgments} 
The authors would like to thank Andrea Aveni, Jordan Brian, Pierpaolo De Blasi, Federico Ferrari, Jürgen Jost, Jeremias Knoblauch, Rostislav Matveev, Vittorio Orlandi, Sonia Petrone, Mai Zhou, Alessandro Zito, and one anonymous reviewer for their helpful comments.  Michele Caprio would like to acknowledge funding from the following grants: NIH 1R01MH118927-01, ARO MURI W911NF2010080. Sayan Mukherjee would like to acknowledge funding from NSF DEB-1840223, NIH R01 DK116187-01, HFSP RGP0051/2017, NSF DMS 17-13012, and NSF CCF-1934964.

\appendix

\section{Further results}\label{app1}
\subsection{Joint distribution of admixture components}
\label{infiniteex}

In this section, we assume that the number $L$ of admixture components in \eqref{fam_first} is known, but the components $\{f_1,\ldots,f_L\}$ are not. We assume they are identically distributed random vectors, but we do not require independence. After realizing that collection $\{f_1,\ldots,f_L\}$ can be seen as a finite exchangeble sequence, we inspect how to approximate its joint distribution applying de Finetti's theorem and a result by Diaconis and Freedman \cite{diaconis1980}. 

Following \cite{aldous_2010}, we can state de Finetti's result from a functional analytic viewpoint in our framework as follows. Let $\mathcal{S} \equiv \Delta^{J-1} \subset \mathbb{R}^{J}$, and recall that a sequence of random variables $X_i$'s is exchangeable if
$$(X_i)_{i \geq 1} \overset{\text{d}}{=} (X_{\text{perm}(i)})_{i \geq 1} \text{,}$$
for any finite permutation $\text{perm}$, where $\overset{\text{d}}{=}$ denotes equality in distribution. We can assume that the elements ${f}_1,\ldots,$ ${f}_{L}$ form a {finite} exchangeable sequence because the order in which they appear provides no additional information about the finite admixture model.

Let $\mathscr{P}(\mathcal{S})\equiv \mathscr{P}(\mathcal{S},\mathcal{B}(\mathcal{S}))$ be the set of probability measures on $(\mathcal{S},\mathcal{B}(\mathcal{S}))$, where $\mathcal{B}(\mathcal{S})$ is the Borel sigma-algebra for $\mathcal{S}$. Let then $\mathscr{P}( \mathscr{P}(\mathcal{S}) )$ be the set of probability measures on $\mathscr{P}(\mathcal{S})$. When we define an infinite exchangeable sequence of $\mathcal{S}$-valued random variables, we are actually defining an exchangeable measure $\Theta\in\mathscr{P}(\mathcal{S}^\infty)$, where $\Theta$ is the distribution of the sequence. Consider the set $\mathfrak{M}:=\{ \mu^\infty:=\mu \times \mu \times \cdots \text{ s.t. } \mu \in \mathscr{P}(\mathcal{S}) \} \subset \mathscr{P}(\mathcal{S}^\infty)$, that is the set of extremal elements of the convex set $\mathfrak{C}$ of exchangeable elements of $\mathscr{P}(\mathcal{S}^\infty)$. Then, we have
\begin{align}\label{defin_funct}
\Theta(A)=\int_{\mathscr{P}(\mathcal{S})} \mu^\infty(A) \text{ } \Lambda(\text{d}\mu), \quad \forall A \subset \mathcal{S}^\infty,
\end{align}
where $\mu^\infty \in \mathfrak{M}$ and we call $\Lambda$ the (unique) \textit{de Finetti measure}. Hence, there is a bijection between $\Lambda \in \mathscr{P} ( \mathscr{P}(\mathcal{S}) )$ and $\Theta \in \mathscr{P}(\mathcal{S}^\infty)$. 

Notice how (functional analytic) de Finetti's theorem \eqref{defin_funct} is similar to Theorem \ref{choq2}. There exists a unique measure ($\Lambda$ for de Finetti and $\nu_c$ for Choquet) supported on the extremal elements ($\mathfrak{M}$ for de Finetti and $E$ for Choquet) of a convex set ($\mathfrak{C}$ for de Finetti and $C$ for Choquet) that allows to represent any element ($\Theta$ for de Finetti and affine function $f(c)$, $c\in C$, for Choquet) within that set.



As we pointed out before, we can assume $\mathbf{f}:=({f}_1,\ldots,$ ${f}_{L})$ to be a {finite} exchangeable sequence. We assumed that the admixture components are identically distributed but not necessarily independent, so let $f_1,\ldots,f_L \sim \mu$. Suppose, without loss of generality, that $\mathbf{f}$ is part of a much longer sequence of $m$ components
$$\left( {f}_1,\ldots,{f}_{L},\ldots,{f}_{m} \right).$$ 
Then, we can use \cite[Theorem 13]{diaconis1980} to compute an approximation of $\Theta_L$, the distribution of $\mathbf{f}$. Let us denote by $\Theta_m$ the distribution of $\left( {f}_1,\ldots,{f}_{L},\ldots,{f}_{m} \right)$; it is an exchangeable probability on $\mathcal{S}^m$. Then, $\Theta_L$, $L \leq m$, is the projection of $\Theta_m$ onto $\mathcal{S}^L$. Define the value $\beta(m,L)$ as
$$\beta(m,L):=1-\frac{m^{-L}m!}{(m-L)!},$$
and notice that $\beta(m,L) \leq \frac{1}{2} \frac{L(L-1)}{m}$.

The theorem states that there exists $\tilde{\Lambda} \in \mathscr{P}(\mathscr{P}(\mathcal{S}))$ such that the probability $\Theta_{\mu L}$ defined on $\mathcal{S}^L$ as
$$\Theta_{\mu L}(A)=\int_{\mathscr{P}(\mathcal{S})} \mu^{L}(A) \text{ } \tilde{\Lambda}(\text{d}\mu) , \quad \forall A \subset \mathcal{S}^L$$
is such that
$d_{TV}(\Theta_L,\Theta_{\mu L}) \leq \beta(m,L)$, for all $L \leq m$. 
We denoted by $\mu^L$ the distribution of $L$ independent picks from $\mu$, that is, $\mu^L((y_1,\ldots,y_L))=\prod_{j=1}^L \mu(y_j)$, and by $d_{TV}$ the total variation distance $$d_{TV}(\Theta_L,\Theta_{\mu L}):=\sup\limits_{A \subset S^L} \left| \Theta_L(A)-\Theta_{\mu L}(A) \right| \text{.}$$

Notice that $\tilde{\Lambda}$ depends on $m$ and $\Theta_m$, but not on $L$, and its analytical form is given in \cite[Proof of Theorem 13]{diaconis1980}.
\subsection{Number of extremal elements of the convex hull having the least amount of vertices}\label{app_2}

The following is an interesting result dealing with the number of extremal elements of a convex hull in $\Delta^{J-1}$ -- but not in any smaller-dimensional unit simplex -- having the least amount of vertices.

\begin{proposition}\label{rema2}
Call $\mathscr{K}\subset \Delta^{J-1}$, $J\in\mathbb{N}$, a polytope such that
\begin{align}\label{ex_within}
\begin{split}
    \tilde{e}:=\#\text{ex}(\mathscr{K}) = \min_{n\in\mathbb{N}} \text{ } &n\\
    \text{subject to} \text{ } &\nexists q\in\{2,\ldots,J\} : \mathscr{K} \subset \Delta^{J-q}
\end{split}
\end{align}
Then, $\tilde{e}=J$.
\end{proposition}

\section{Proofs}\label{proofs}

Recall that a chain $\mathcal{F}_0(P)\subset \mathcal{F}_1(P) \subset \cdots \subset \mathcal{F}_{h}(P)$ of $i$-dimensional faces of a polytope $P$ is called a \textit{tower} (or a \textit{flag}) of $P$, and we denote it as $\tau(P)$.

\begin{proof}[Proof of Theorem \ref{growth}]
 In \cite[Theorem 6]{reitzner} and \cite[Theorem 5]{barany2}, the authors show that, given a convex polytope $P$ in $\mathbb{R}^d$, if we call $P_n$ the convex hull of $n$ points sampled iid from a uniform on $P$, then 
$$\mathbb{E} \left[F_0(P_n)\right] = \frac{\tau(P)}{(d+1)^{d-1}  (d-1)!}  (\log n)^{d-1} + \mathcal{O} \left( (\log n)^{(d-2)} \log\log n \right).$$
Then, since $\Delta^{J-1}$ is a convex polytope in $\mathbb{R}^J$, and given the way we defined $K_n$, equations \eqref{first_theorem} and \eqref{equation_imp_first} follow immediately after realizing that $\tau(\Delta^{J-1})=(J-1)! \cdot J$.
\end{proof}

\begin{proof}[Proof of Theorem \ref{var_thm}]

Let $\text{Vol}$ denote the volume operator. In \cite[Theorem 1.3]{reitzner_var}, the authors show that, given a convex polytope $P$ in $\mathbb{R}^d$ such that $\text{Vol}(P)=1$, if we call $P_n$ the convex hull of $n$ points sampled iid from a uniform on $P$, then there exists a constant $C>0$ such that
$$C \tau(P)(\log n)^{d-1} < \mathbb{V}[F_0(P_n)] < C \tau(P)^3(\log n)^{d-1}.$$
Recall that $\text{Vol}(\Delta^{J-1})=\sqrt{J}/[(J-1)!]$. Then, since $\text{Vol}(\Delta^{J-1})$ is in a fixed relation 
$$x \mapsto f(x)=\frac{\sqrt{J}}{(J-1)!} \cdot x$$ 
with $\text{Vol}(P)$, because $\Delta^{J-1}$ is a convex polytope in $\mathbb{R}^J$, and given the way we defined $K_n$, equation \eqref{var_eq} follows immediately.
\end{proof}

\begin{proof}[Proof of Theorem \ref{clt_2}]
In \cite[Corollary 1.2]{pardon}, the author shows that, given a convex polytope $P$ in $\mathbb{R}^2$ of unit area, if we call $P_n$ the convex hull of $n$ points sampled iid from a uniform on $P$, then
$$\lim_{n\rightarrow\infty} \sup_{x\in\mathbb{R}} \left| \mathbb{P}\left( \frac{F_0(P_n)-\mathbb{E}[F_0(P_n)]}{\sqrt{\mathbb{V}[F_0(P_n)}} \leq x \right) - \Phi(x) \right|=0.$$
Recall that the area of the unit simplex in $\mathbb{R}^3$ is $\sqrt{3}/2$. Then, since the area of $\Delta^2$ is in a fixed relation $x\mapsto f(x)=\frac{\sqrt{3}}{2}x$ with the area of $P$, because $\Delta^2$ is a convex polytope in $\mathbb{R}^2$, and given the way we defined $K_n$, equation \eqref{clt_2_eq} follows immediately.
\end{proof}

\begin{proof}[Proof of Theorem \ref{clt_gen}]
In \cite[Theorems 2, 6]{reitzner1} the author shows that, given a smooth compact convex set $P$ in $\mathbb{R}^d$, if we call $P_n$ the convex hull of $n$ points sampled iid from a uniform on $P$, then
$$\left| \mathbb{P} \left( F_0(P_n)\leq \mathbb{E}[F_0(P_n)]+x\sqrt{\mathbb{V}[F_0(P_n)]} \right) - \Phi(x) \right|= \mathcal{O}\left( n^{-\frac{1}{2(J+1)}}\left(\log n\right)^{2+\frac{2}{J+1}} \right).$$
Since $\hat{\Delta}^{J-1}_\epsilon$ is a smooth compact convex set in $\mathbb{R}^J$, and given the way we defined $\hat{K}_n$, equation \eqref{clt_gen_eq} follows immediately.
\end{proof}

\begin{proof}[Proof of Theorem \ref{conc_unif_thm}]
    In \cite[Theorem 2.11, Section 7]{vu}, the author shows that given a smooth compact convex set $P$ in $\mathbb{R}^d$, if we call $P_n$ the convex hull of $n$ points sampled iid from a uniform on $P$, then there exist positive constants $c,\Xi,\epsilon_0$ such that for any $\epsilon\in(0,\epsilon_0]$, $V\geq \Xi n^{\frac{d-1}{d+1}}$, $C\geq n\epsilon$, and $\lambda\in(0,\frac{V}{4C^2})$, the following holds
        $$\mathbb{P}\left( \left| F_0(P_n) -\mathbb{E}[F_0(P_n)] \right| \geq \sqrt{\lambda V} \right) \leq 2 \exp(-\lambda/4)+\exp(-c\epsilon n) + \exp\left(-cn^\frac{d-1}{3d+5}\right).$$
        Since $\hat{\Delta}^{J-1}_\epsilon$ is a smooth compact convex set in $\mathbb{R}^J$, and given the way we defined $\hat{K}_n$, equation \eqref{conc_unif_eq} follows immediately.
\end{proof}

\begin{proof}[Proof of Theorem \ref{extrema2}]
By hypothesis, we have that for all $n\geq N$, $\mathbb{E}[T(n)]=\gamma_n \mathbb{E}[M(n)]$.
In addition, 
by Theorem \ref{growth} we have that 
$$\lim_{n \rightarrow \infty} \frac{\mathbb{E}[M(n)]}{(\log n)^{J-1}} =c(J).$$
Hence we obtain that
$$\lim_{n\rightarrow\infty} \frac{\mathbb{E}[T(n)]}{{\gamma_n}(\log n)^{J-1}}=\lim_{n\rightarrow\infty} \frac{\gamma_n\mathbb{E}[M(n)]}{{\gamma_n}(\log n)^{J-1}}= \lim_{n\rightarrow\infty} \frac{\mathbb{E}[M(n)]}{(\log n)^{J-1}}= c(J),$$
concluding the proof.
\end{proof}

\begin{proof}[Proof of Corollary \ref{extrema2_cor}]
Suppose that the assumptions of Theorem \ref{extrema2} hold and that $\gamma_n=\mathcal{O}(\varpi_n)$. This latter means that there exists $R\in\mathbb{R}$ and $N\in\mathbb{N}$ such that for all $n\geq N$,
$$\frac{\gamma_n}{\varpi_n} \leq R.$$
Then,
$$\lim_{n\rightarrow\infty} \frac{\mathbb{E}[T(n)]}{\varpi_n (\log n)^{J-1}}=\lim_{n\rightarrow\infty} \frac{\gamma_n \mathbb{E}[M(n)]}{\varpi_n (\log n)^{J-1}}=\lim_{n\rightarrow\infty} \frac{\gamma_n}{\varpi_n} \lim_{n\rightarrow\infty} \frac{\mathbb{E}[M(n)]}{(\log n)^{J-1}} \leq R c(J),$$
concluding the proof. 

\end{proof}

\begin{proof}[Proof of Proposition \ref{simplex_conseq}]
The proposition is an immediate consequence of 
Theorem \ref{choq2}.
\end{proof}


\begin{proof}[Proof of Theorem \ref{thm1}]
The uniform almost sure convergence statement is a consequence of the Glivenko-Cantelli's Theorem \cite{vdv_asy}, while the rate of convergence can be easily derived from \cite[Lemma 2.1]{zhou}, via a Borel-Cantelli argument.\footnote{We refer the reader that is further interested in the rates of the Glivenko-Cantelli convergence to the recent work \cite{dolera}.}
\end{proof}
%

\begin{proof}[Proof of Theorem \ref{algo_cons}]
    First notice that $L_T$ depends on the amount $n$ of data available to perform the estimating procedure in Algorithm \ref{alg:2step}, so we can write $L_T\equiv L_T(n)$. Suppose now for the sake of contradiction that $\tilde{M}:=\lim_{n\rightarrow\infty} L_T(n) \neq M^\star$. Then, we have two cases, either $\tilde{M} > M^\star$, or $\tilde{M} < M^\star$.

    Case 1: If $\tilde{M} > M^\star$, then there exist a collection $\{\hat{f}_1,\ldots,\hat{f}_{\tilde{M}-M^\star}\}$ of admixture components that can be written as a convex combination of the remaining $M^\star$ components. Together with the assumption of the MLE in \eqref{mixtureml} being consistent, this contradicts Step 2 of Algorithm \ref{alg:2step}.

    Case 2: Suppose now $\tilde{M} < M^\star$. Then, since by Algorithm \ref{alg:2step} $L_T(n) \leq L_0$, for all $n$, we have that $\tilde{M}\leq L_0$. So, if $\tilde{M} < M^\star$, it follows that either $L_0 < M^\star$, or $L_0 \geq M^\star$. If $L_0 \geq M^\star$, Algorithm \ref{alg:2step} would have stopped at $\tilde{M} = M^\star$. But since we are assuming $\tilde{M} < M^\star$, then this means that $L_0 < M^\star$, which contradicts the assumption of our theorem. This concludes the proof.
\end{proof}

\begin{proof}[Proof of Proposition \ref{rema2}]
Suppose for the sake of contradiction that $\tilde{e}\neq J$. This means that either $\tilde{e}>J$, or $\tilde{e}<J$. If the latter holds, then there exists $q^\prime\in\{2,\ldots,J\}$ such that $\mathscr{K} \subset \Delta^{J-q^\prime}$, which contradicts \eqref{ex_within}. If instead $\tilde{e}>J$, then we can find $\mathscr{K}^\prime \subsetneq \mathscr{K}$ such that $\# \text{ex}(\mathscr{K}^\prime)<\tilde{e}$, but $\mathscr{K}^\prime$ is still a proper subset of $\Delta^{J-1}$, thus again contradicting \eqref{ex_within}. This concludes the proof.
\end{proof}

\bibliographystyle{plain}
\bibliography{References}
\end{document}